\documentclass[review]{elsarticle}
\usepackage[top=3cm,bottom=3cm,left=3cm,right=3cm]{geometry}
\usepackage{amsmath,amssymb,amsthm}
\usepackage{xcolor}
\usepackage{tikz}

\usepackage{array,longtable}
\usepackage{algorithm}
\usepackage{algpseudocode}
\biboptions{numbers,sort&compress}
\usepackage{hyperref}
\definecolor{edgezero}{HTML}{2166AC}
\definecolor{edgeone}{HTML}{D66000}
\definecolor{edgetwo}{HTML}{008060}
\tikzset{av/.style={circle,draw,fill=white,inner sep=2pt,minimum size=7mm},
 bv/.style={rectangle,draw,fill=white,inner sep=2pt,minimum size=7mm},
 ee/.style={line width=1.2pt}, every node/.style={font=\small}}
\newtheorem{theorem}{Theorem}
\newtheorem{lemma}[theorem]{Lemma}
\newtheorem{proposition}[theorem]{Proposition}
\newtheorem{corollary}[theorem]{Corollary}
\newtheorem{conjecture}[theorem]{Conjecture}
\theoremstyle{definition}

\newtheorem{problem}[theorem]{Problem}
\newtheorem{claim}[theorem]{Claim}
\theoremstyle{remark}

\tikzset{
  graph/.style={
    every node/.style={font=\small,inner sep=1pt,text=black},
    every label/.style={font=\small,inner sep=1pt,text=black},
    label distance=2pt,
    line width=.55pt
  },
  graphvertex/.style={
    circle,draw=black,fill=white,minimum size=4.5pt,
    inner sep=0pt,line width=.5pt
  },
  graphselected/.style={graphvertex,fill=black},
  graphcovered/.style={graphvertex,fill=gray!50},
  graphuncovered/.style={graphvertex},
  graphunknown/.style={graphvertex,fill=gray!20,densely dashed},
  residualedge/.style={draw=black,line width=1.2pt},
  graphguide/.style={draw=gray,densely dotted,rounded corners=3pt,line width=.45pt},
  pics/isolationlegend/.style={code={
    \node[graphselected] at (0,0) {};
    \node[anchor=west] at (.22,0) {selected};
    \node[graphcovered] at (2.35,0) {};
    \node[anchor=west] at (2.57,0) {covered};
    \node[graphuncovered] at (4.65,0) {};
    \node[anchor=west] at (4.87,0) {uncovered};
  }}
}

\begin{document}
\begin{frontmatter}
\title{$K_{1,2}$-Isolation in Triangle-Free Cubic Graphs\tnoteref{funding}}
\tnotetext[funding]{This work is supported by the Scientific and Technological Research Council of Turkey (TUBITAK) under grant no.\ 126F284.}
\author[gtu]{Yasemin B\"uy\"uk\c{c}olak}
\ead{y.buyukcolak@gtu.edu.tr}
\author[gtu]{Fatih Yetgin}
\ead{fyetgin@gtu.edu.tr}
\address[gtu]{Department of Mathematics, Gebze Technical University, Kocaeli, Turkey}

\begin{abstract}
A $K_{1,2}$-isolating set of a graph $G$ is a vertex set
whose closed neighborhood can be removed to leave
a graph of maximum degree at most one.
The minimum cardinality of such a set is denoted by $\iota_1(G)$.
We establish the sharp bounds
$\lceil \frac{n}{7}\rceil\leq\iota_1(G)\leq\lfloor \frac{n}{4}\rfloor$
for connected triangle-free cubic graphs of order $n$,
with both bounds attained.
The upper bound confirms the triangle-free case of
a conjecture for cubic graphs.
The lower bound holds more generally for connected cubic graphs,
and we characterize those satisfying $\iota_1(G)=\frac{n}{7}$, including
the triangle-free members.
We further characterize all connected bipartite cubic graphs
satisfying $\iota_1(G)=\frac{n}{4}$, showing that they are precisely
the members of an infinite family together with the
M\"obius-Kantor graph.
We also prove that these graphs are vertex-critical:
deleting any vertex decreases their
$K_{1,2}$-isolation number by exactly one.
In addition, we construct a nonbipartite
triangle-free cubic graph of order sixteen and girth five
attaining the upper bound.
\end{abstract}

\end{frontmatter}

\section{Introduction}
Throughout the paper, all graphs are finite, simple, and undirected.
For a graph $G$, let $V(G)$ and $E(G)$ denote
its vertex and edge sets, respectively, and let
$n=|V(G)|$ be its order.
We write $d_G(v)$ for the degree of a vertex
$v$, and $\delta(G)$ and $\Delta(G)$ for the
minimum and maximum degrees of $G$.
For $v\in V(G)$, let $N_G(v)$ and
$N_G[v]=N_G(v)\cup\{v\}$ denote its open
and closed neighborhoods, respectively.
For $S\subseteq V(G)$, put
$N_G[S]=\bigcup_{v\in S}N_G[v]$.
Subscripts are omitted when the graph is clear from the context.
We write $G[X]$ for the subgraph induced
by $X\subseteq V(G)$ and
$G-X=G[V(G)\setminus X]$.
The path and cycle on $k$ vertices are
denoted by $P_k$ and $C_k$, respectively.
A graph is \emph{cubic} if every vertex
has degree three and \emph{triangle-free}
if it contains no $C_3$.
For disjoint vertex sets
$X,Y\subseteq V(G)$, let $e_G(X,Y)$
denote the number of edges joining
$X$ and $Y$.
A partition $V(G)=A\sqcup B$ determines
a \emph{cut}, consisting of all edges
between $A$ and $B$.
A \emph{maximum cut} is one containing
the largest possible number of edges.
A \emph{matching} is a set of pairwise
vertex-disjoint edges.

Graph isolation generalizes domination by restricting
the subgraphs that may remain outside the closed
neighborhood of a selected vertex set.  In the framework introduced by Caro and
Hansberg~\cite{CaroHansberg2017}, a set $S\subseteq V(G)$ is
\emph{$\mathcal F$-isolating} if $G-N[S]$ contains no member of
the graph family $\mathcal F$ as a subgraph. The minimum cardinality of such a set is
the \emph{$\mathcal F$-isolation number},
denoted by $\iota(G,\mathcal F)$.
When $\mathcal F=\{H\}$, we write
$\iota_H(G)$ instead of $\iota(G,\mathcal F)$.
In particular, when $H=K_{1,k+1}$,
we use the notation $\iota_k(G)$ and
refer to this parameter as the
\emph{$k$-isolation number}.
Several familiar parameters arise as
special cases of this framework.
For $\mathcal F=\{K_1\}$, an
$\mathcal F$-isolating set is precisely
a \emph{dominating set}, that is, a set
$S$ satisfying $N_G[S]=V(G)$.
Its minimum cardinality is the
\emph{domination number} $\gamma(G)$.
For $\mathcal F=\{K_2\}$, we obtain
\emph{ordinary isolation}, which coincides
with \emph{vertex-edge domination}: a set $S$
is isolating if $G-N_G[S]$ is edgeless.
The minimum cardinality of such a set
is the \emph{isolation number}
$\iota(G)=\iota_0(G)$.

\medskip

In this paper, we study $K_{1,2}$-isolation where $K_{1,2}\cong P_3$,
also known as $P_3$-isolation \cite{WeiZhangZhao2024} or $3$-path isolation
\cite{BartoloBorgScicluna2026,
BartoloBorgScicluna2026Paths}. A set $S$ is \emph{$K_{1,2}$-isolating}
if $G-N_G[S]$ contains no $P_3$ as
a subgraph, or equivalently, has
maximum degree at most one.
The minimum cardinality of such a set
is the \emph{$K_{1,2}$-isolation number}
$\iota_1(G)$.
Unlike ordinary isolation, which requires
the remaining vertices to be independent,
$K_{1,2}$-isolation allows isolated edges
to remain. An isolating set is \emph{independent}
if no two of its vertices are adjacent,
and $\iota^{\mathrm i}(G)$ denotes
the minimum cardinality of an
independent ordinary isolating set.
A vertex is said to be \emph{covered}
by $S$ if it belongs to $N_G[S]$;
otherwise, it is \emph{uncovered}.
These definitions imply
$\iota_1(G)\leq\iota(G)\leq\gamma(G)$
and $\iota(G)\leq\iota^{\mathrm i}(G)$. All of these inequalities hold with equality even for cubic graphs,
as shown by $K_4$, for which
$
\iota_1(K_4)=\iota(K_4)=\gamma(K_4)=\iota^{\mathrm i}(K_4)=1
$.

\medskip

The following result establishes the additivity of $\mathcal F$-isolation over connected components.

\begin{lemma}[{\cite{ZhangWu2021}}]\label{lem:component-additivity}
Let $\mathcal{F}$ be a family of connected graphs, and let
$G_1,\ldots,G_s$ be the distinct connected components of a graph $G$.
Then
$
\iota(G,\mathcal{F})
=\sum_{j=1}^{s}\iota(G_j,\mathcal{F})
$.
\end{lemma}
\noindent
Lemma~\ref{lem:component-additivity}
allows us to restrict attention to connected graphs.
Hence, the graph $G$ in our results is assumed to be connected;
auxiliary graphs and subgraphs need not be connected.
The corresponding statements for disconnected graphs follow by applying
each result to its connected components, under the same hypotheses.

\medskip

A natural question in graph isolation is to determine upper bounds
on the isolation number in terms of the order of a graph.
For $K_{1,2}$-isolation, Zhang and Wu~\cite{ZhangWu2021}
established the following sharp bounds.

\begin{theorem}[\cite{ZhangWu2021}]
Let $G$ be a connected graph of order $n$. Then:
\begin{enumerate}
    \item If $G \notin \{P_3,C_3,C_6\}$, then
    $\iota_1(G) \leq \frac{2n}{7}$.

    \item If $G$ has girth at least seven and
    $G \notin \{P_3,C_7,C_{11}\}$, then
    $\iota_1(G) \leq \frac{n}{4}$.
\end{enumerate}
\end{theorem}
These results naturally raise the question of whether the
stronger bound $\frac{n}{4}$ holds for other graph classes
without imposing a girth condition.
In particular, Yin, An and Wu~\cite{YinAnWu2024}
established this bound for claw-free cubic graphs
and conjectured that it holds for all cubic graphs.
\begin{conjecture}[\cite{YinAnWu2024}]
\label{conj:cubic-path-isolation}
Every cubic graph $G$ of order $n$ satisfies
$\iota_1(G) \leq \frac{n}{4}$.
\end{conjecture}

Related results establish the bound in
Conjecture~\ref{conj:cubic-path-isolation}
under various restrictions on short cycles.
Wei, Zhang and Zhao~\cite{WeiZhangZhao2024}
proved that every connected triangle-free graph
without induced six-cycles satisfies
$\iota_1(G)\leq \lfloor \frac{n}{4}\rfloor$,
except for $P_3$, $C_7$, and $C_{11}$.
Huang, Zhang and Jin~\cite{HuangZhangJin2026}
established the same bound for connected graphs
containing no six-cycle, or containing neither
an induced five-cycle nor an induced six-cycle,
with the exceptions $P_3$, $C_3$, $C_7$,
and $C_{11}$.
More recently, Bartolo, Borg and
Scicluna~\cite{BartoloBorgScicluna2026}
proved that connected subcubic graphs without
induced six-cycles satisfy
$\iota_1(G)\leq \lfloor \frac{n}{4}\rfloor$,
apart from twelve exceptional graphs,
none of which is cubic.
In a subsequent paper~\cite{BartoloBorgScicluna2026Paths},
they obtained the bound
$\iota_1(G)\leq \lfloor\frac{n+1}{4}\rfloor$
for all connected graphs without induced six-cycles.
Since every cubic graph has even order,
this confirms the conjecture for cubic graphs
without induced six-cycles.
The conjecture has also been confirmed for
bipartite cubic graphs, without any restriction
on six-cycles.
Boyer and Goddard~\cite{BoyerGoddard2025}
showed that every bipartite cubic graph satisfies
$\iota^{\mathrm i}(G)\leq \frac{n}{4}$.
Together with
$\iota_1(G)\leq\iota(G)\leq\iota^{\mathrm i}(G)$,
this establishes the conjectured bound for
this class as well.

The preceding results leave open the conjecture for connected
nonbipartite triangle-free cubic graphs containing induced six-cycles.
In this paper, we resolve this remaining case by proving
that the conjectured upper bound holds for all connected triangle-free
cubic graphs.
In addition, we give a lower bound that holds for
all connected cubic graphs and remains sharp even within the
triangle-free class.
Together, we prove the following sharp bounds as
the main result of this paper.

\begin{theorem}\label{thm:main-bounds}
Let $G$ be a connected triangle-free cubic graph of order $n$. Then
\[
\left\lceil \frac{n}{7} \right\rceil
\leq \iota_1(G)
\leq \left\lfloor \frac{n}{4} \right\rfloor.
\]
Furthermore, the lower bound holds for every connected cubic graph.
\end{theorem}
We also investigate equality in the linear bounds $\frac{n}{7}$ and $\frac{n}{4}$.
For the lower bound, we provide an explicit construction of all
connected cubic graphs satisfying $\iota_1(G)=\frac{n}{7}$
and identify the triangle-free members of this family.
For the upper bound, we completely characterize the
connected bipartite cubic graphs satisfying
$\iota_1(G)=\frac{n}{4}$, showing that they are precisely
the members of an infinite family $\mathcal H$
and the M\"obius-Kantor graph.
Moreover, we prove that these graphs are vertex-critical:
deleting any vertex from such a graph decreases
its $K_{1,2}$-isolation number by one.
We also construct a nonbipartite triangle-free cubic graph
of order sixteen and girth five attaining the upper bound,
demonstrating that equality is not confined to
bipartite graphs.

\medskip

The rest of the paper is organized as follows.
In Section~\ref{sec:triangle-free-path-isolation},
we establish the lower and upper bounds on the
$K_{1,2}$-isolation number.
Section~\ref{sec:equality} investigates the equality cases,
characterizing connected cubic graphs with $\iota_1(G)=\frac{n}{7}$
and connected bipartite cubic graphs with $\iota_1(G)=\frac{n}{4}$.
We also establish the vertex-criticality of the latter graphs.
Finally, Section~\ref{sec:conclusion} presents a
nonbipartite graph attaining the upper bound
and discusses some open problems.

\section{Linear Bounds on the \texorpdfstring{$K_{1,2}$}{K1,2}-Isolation Number}
\label{sec:triangle-free-path-isolation}

In this section, we prove the bounds stated in
Theorem~\ref{thm:main-bounds}.
The lower bound follows from a counting argument
valid for all connected cubic graphs.
For the upper bound, we use a maximum cut to obtain
a spanning bipartite subgraph with suitable structural
properties, allowing us to apply a result on total domination.

\medskip

We first record elementary consequences of cut maximality.

\begin{lemma}\label{lem:triangle-free-maximum-cut}
Let $G$ be a connected triangle-free cubic graph,
and let $V(G)=A\sqcup B$ define a maximum cut.
Let $F$ be the spanning subgraph consisting of the
cut edges, and let $R$ be the spanning subgraph with
$E(R)=E(G)\setminus E(F)$.
Then $F$ is connected and bipartite, with
$2\leq d_F(v)\leq 3$ for every vertex and
$\Delta(F)=3$.
Moreover, $R$ is a matching together with isolated
vertices, and every component of $F[L]$, where
$L=\{v\in V(G):d_F(v)=2\}$, is a $P_1$ or a $P_2$.
\end{lemma}

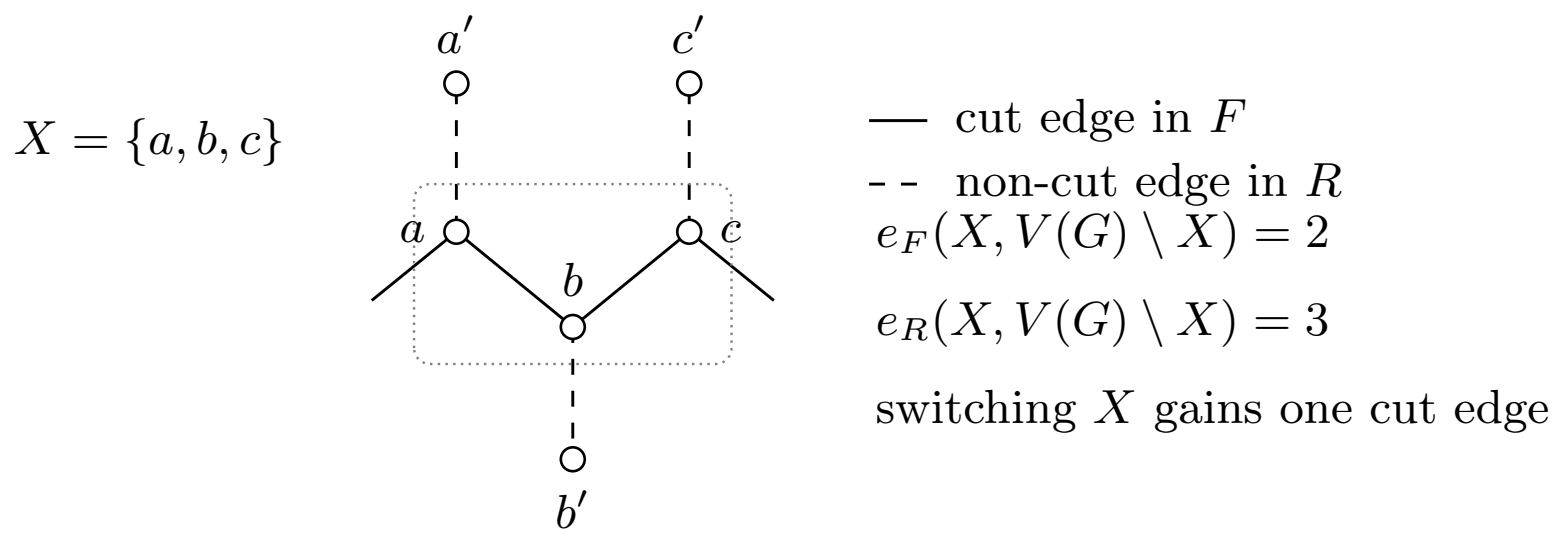
\begin{figure}[tbp]
\centering
\begin{tikzpicture}[graph,scale=0.7]
\coordinate (a) at (-1.1,.25);
\coordinate (b) at (0,-.65);
\coordinate (c) at (1.1,.25);
\coordinate (am) at (-1.1,1.65);
\coordinate (bm) at (0,-1.9);
\coordinate (cm) at (1.1,1.65);
\draw (a)--(b)--(c);
\draw (a)--(-1.9,-.4) (c)--(1.9,-.4);
\draw[dashed] (a)--(am) (b)--(bm) (c)--(cm);
\draw[graphguide] (-1.5,-1) rectangle (1.5,.7);
\node at (-4,1.1) {$X=\{a,b,c\}$};
\node[graphvertex,label=left:$a$] at (a) {};
\node[graphvertex,label=above:$b$] at (b) {};
\node[graphvertex,label=right:$c$] at (c) {};
\node[graphvertex,label=above:$a'$] at (am) {};
\node[graphvertex,label=below:$b'$] at (bm) {};
\node[graphvertex,label=above:$c'$] at (cm) {};
\draw (2.8,1.3)--(3.35,1.3);
\node[anchor=west] at (3.55,1.3) {cut edge in $F$};
\draw[dashed] (2.8,.7)--(3.35,.7);
\node[anchor=west] at (3.55,.7) {non-cut edge in $R$};
\node[anchor=west,align=left] at (2.8,-.6)
  {$e_F(X,V(G)\setminus X)=2$\\
   $e_R(X,V(G)\setminus X)=3$\\
   switching $X$ gains one cut edge};
\end{tikzpicture}
\caption{The three-vertex cut exchange}
\label{fig:triangle-free-cut-exchange}
\end{figure}

\begin{proof}
For any $X\subseteq V(G)$, consider switching all vertices
of $X$ to the opposite side of the cut.
Only edges with exactly one endpoint in $X$ change
their status: cut edges become non-cut edges, and
non-cut edges become cut edges.
Since the original cut is maximum, this switch
cannot increase its size. Hence
\begin{equation}\label{eq:maximum-cut-exchange}
e_F(X,V(G)\setminus X)
\geq e_R(X,V(G)\setminus X)
\qquad\text{for every }X\subseteq V(G).
\end{equation}
Taking $X=\{v\}$ gives $d_F(v)\geq d_R(v)$.
Since $d_F(v)+d_R(v)=3$, we have
$d_F(v)\geq 2$ and $d_R(v)\leq 1$.
Thus $R$ is a matching together with isolated vertices.

We next show that $F$ is connected.
Otherwise, the vertex set $X$ of a component of $F$
would satisfy $e_F(X,V(G)\setminus X)=0$.
Since $G$ is connected, at least one edge of $R$
must join $X$ to its complement, contradicting
\eqref{eq:maximum-cut-exchange}.
Therefore, $F$ is connected.
It is bipartite by construction.

Let $L=\{v\in V(G):d_F(v)=2\}$.
Suppose that $abc$ is a path in $F[L]$, and put $X=\{a,b,c\}$.
Triangle-freeness gives $E(G[X])=\{ab,bc\}$.
Since $d_F(v)=2$ and $d_R(v)=1$ for every $v\in X$,
exactly two edges of $F$ and three edges of $R$ leave $X$,
contradicting \eqref{eq:maximum-cut-exchange}.
Figure~\ref{fig:triangle-free-cut-exchange} illustrates this exchange.

It follows that $F[L]$ contains no path
on three vertices, so every component of
$F[L]$ is a $P_1$ or a $P_2$.
Finally, if $\Delta(F)\leq 2$, then
$d_F(v)=2$ for every vertex, giving
$L=V(G)$ and $F[L]=F$.
This contradicts the structure of $F[L]$
established above.
Hence $\Delta(F)=3$.
\end{proof}

The preceding lemma provides the structural conditions
needed to apply a known bound on total domination.
Recall that a set $D\subseteq V(G)$ is
\emph{total dominating} if every vertex of $G$
has a neighbor in $D$.
For a graph without isolated vertices, the minimum
cardinality of such a set is the
\emph{total domination number}, denoted by
$\gamma_{\mathrm t}(G)$.
We use the following result of Lam and
Wei~\cite{LamWei2007}, also stated by
Henning and Yeo~\cite{HenningYeo2007}.

\begin{proposition}[\cite{LamWei2007, HenningYeo2007}]
\label{prop:td-short-degree-two-components}
Let $H$ be a connected graph of order $m$
with $\delta(H)\geq 2$.
If every component of the subgraph induced
by the degree-two vertices of $H$ has
order at most two, then
$\gamma_{\mathrm t}(H)\leq m/2$.
\end{proposition}

By Lemma~\ref{lem:triangle-free-maximum-cut},
the spanning bipartite subgraph $F$ satisfies
the hypotheses of Proposition \ref{prop:td-short-degree-two-components}.
We can therefore use a small total dominating
set of $F$ to construct a $K_{1,2}$-isolating
set of $G$.
We now prove the main theorem.

\begin{proof}[Proof of Theorem~\ref{thm:main-bounds}]
We first establish the lower bound for every connected cubic graph,
without assuming triangle-freeness.
Let $G$ be such a graph of order $n$, and let
$S$ be a minimum $K_{1,2}$-isolating set.
Put $T=N_G[S]\setminus S$ and
$U=V(G)\setminus N_G[S]$.

Since $G$ is cubic, $|T|\leq 3|S|$.
Every vertex of $U$ has no neighbor in $S$
and at most one neighbor in $U$.
Thus each vertex of $U$ has at least two
neighbors in $T$.
On the other hand, every vertex of $T$
has a neighbor in $S$ and therefore
at most two neighbors in $U$.
Consequently,
\(
2|U|\leq e_G(T,U)\leq 2|T|,
\)
which gives $|U|\leq |T|$.
It follows that
\[
n=|S|+|T|+|U|
\leq |S|+2|T|
\leq 7|S|.
\]
Since $|S|=\iota_1(G)$ is an integer,
we obtain $\iota_1(G)\geq\lceil \frac{n}{7}\rceil$.

For the upper bound, let $G$ be connected and triangle-free,
and choose a maximum cut $V(G)=A\sqcup B$.
Let $F$ and $R$ be the corresponding
spanning subgraphs defined in
Lemma~\ref{lem:triangle-free-maximum-cut}.
By this lemma, $F$ is connected,
$\delta(F)\geq 2$, and every component
induced by its degree-two vertices
has order at most two.
Proposition~\ref{prop:td-short-degree-two-components}
therefore guarantees a total dominating
set $D$ of $F$ with $|D|\leq \frac{n}{2}$.

Interchanging $A$ and $B$ if necessary,
we may assume that $|D\cap A|\leq |D\cap B|$.
Set $S=D\cap A$.
Since $F$ is bipartite and $D$ totally
dominates $F$, every vertex of $B$
has a neighbor in $S$.
Hence $B\subseteq N_G[S]$, so every
vertex remaining in $G-N_G[S]$ belongs to $A$.
Every edge of $G[A]$ belongs to $R$,
which has maximum degree at most one.
Thus $\Delta(G-N_G[S])\leq 1$,
showing that $S$ is a $K_{1,2}$-isolating
set of $G$.
Moreover,
\[
\iota_1(G)\leq |S|
\leq \left\lfloor\frac{|D|}{2}\right\rfloor
\leq \left\lfloor\frac{n}{4}\right\rfloor.
\]

This completes the proof.
\end{proof}

\subsection{The \texorpdfstring{$K_{1,2}$}{K1,2}-Isolation Number of M\"obius Ladders}

The bounds in Theorem~\ref{thm:main-bounds} apply to all
connected triangle-free cubic graphs. For certain graph families,
however, the $K_{1,2}$-isolation number can be determined
exactly. We illustrate this with M\"obius ladders, a
well-known family of connected triangle-free cubic graphs.

For every even integer $n\geq 6$, let $M_n$ be the graph
obtained from the cycle $v_0v_1\cdots v_{n-1}v_0$
by adding the edges $v_iv_{i+\frac{n}{2}}$ for
$0\leq i\leq \frac{n}{2}-1$.
The following result determines its $K_{1,2}$-isolation
number.

\begin{proposition}\label{prop:mobius-isolation}
For every even integer $n\geq 6$, the M\"obius ladder
$M_n$ satisfies
$\iota_1(M_n)=\lceil \frac{n}{6}\rceil$.
\end{proposition}

\begin{proof}
Write $n=2m$, where $m=3q+r$ and $r\in\{0,1,2\}$.
Let $A=\{v_0,\ldots,v_{m-1}\}$ and
$B=\{v_m,\ldots,v_{n-1}\}$.
Each half induces a path, and the edges $v_iv_{i+m}$
join corresponding vertices of the two halves.
The remaining edges between them are
$v_{m-1}v_m$ and $v_{n-1}v_0$.

We first construct a $K_{1,2}$-isolating set
of size $\lceil \frac{n}{6}\rceil$.
Define
\[
S=
\begin{cases}
\{v_{3j}:0\leq j\leq q-1\}, & r=0,\\
\{v_0\}\cup\{v_{3j+1}:0\leq j\leq q-1\}, & r=1,\\
\{v_0\}\cup\{v_{3j+2}:0\leq j\leq q-1\}, & r=2.
\end{cases}
\]
Then $S\subseteq A$, $|S|=\lceil \frac{m}{3}\rceil$,
and $v_0\in S$. Moreover, consecutive selected
indices differ by at most three, and the last
selected vertex is $v_{m-3}$.
Every vertex of $A$, except possibly $v_{m-1}$,
belongs to $N[S]$. If $v_{m-1}$ remains uncovered,
it is isolated in $M_n-N[S]$, since its three
neighbors are covered.
In $B$, the selected vertices cover their
corresponding vertices through the added edges.
Together with the coverage of $v_m$ and $v_{n-1}$
by $v_0$, the gaps between selected indices
ensure that the uncovered vertices of $B$
induce only isolated vertices and edges.
Thus $\Delta(M_n-N[S])\leq 1$, proving that
\(
\iota_1(M_n)\leq |S|
=\left\lceil\frac{n}{6}\right\rceil.
\)
Figure~\ref{fig:mobius-16} illustrates this construction.

\begin{figure}[htbp]
\centering
\begin{tikzpicture}[graph,scale=0.7]
\foreach \i in {0,...,15}{
  \pgfmathsetmacro{\ang}{78.75-22.5*\i}
  \coordinate (v\i) at (\ang:2.25);
  \coordinate (l\i) at (\ang:2.72);
}

\draw[dashed,gray] (0,-2.95)--(0,2.95);
\node at (1.55,3.25) {first half};
\node at (-1.55,3.25) {second half};

\foreach \i in {0,...,15}{
  \pgfmathtruncatemacro{\j}{mod(\i+1,16)}
  \draw (v\i)--(v\j);
}
\foreach \i in {0,...,7}{
  \pgfmathtruncatemacro{\j}{\i+8}
  \draw (v\i)--(v\j);
}

\draw[residualedge] (v11)--(v12);

\foreach \i in {0,2,5}
  \node[graphselected] at (v\i) {};

\foreach \i in {1,3,4,6,8,10,13,15}
  \node[graphcovered] at (v\i) {};

\foreach \i in {7,9,11,12,14}
  \node[graphuncovered] at (v\i) {};

\foreach \i in {0,...,15}
  \node at (l\i) {$v_{\i}$};

\node[graphselected] at (4.2,1.2) {};
\node[anchor=west] at (4.55,1.2) {selected};

\node[graphcovered] at (4.2,0.4) {};
\node[anchor=west] at (4.55,0.4) {covered};

\node[graphuncovered] at (4.2,-0.4) {};
\node[anchor=west] at (4.55,-0.4) {uncovered};

\draw[residualedge] (3.95,-1.2)--(4.45,-1.2);
\node[anchor=west] at (4.55,-1.2) {uncovered edge};

\end{tikzpicture}
\caption{$M_{16}$ with $S=\{v_0,v_2,v_5\}$.}
\label{fig:mobius-16}
\end{figure}
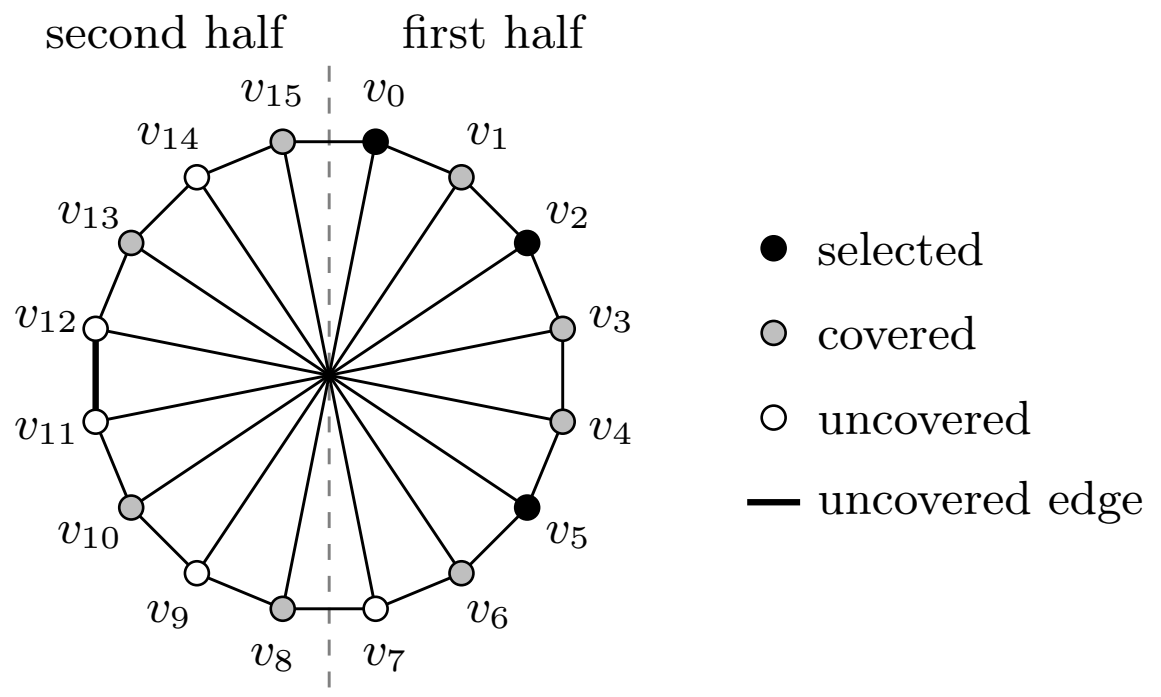

For the reverse inequality, let $S$ be any
$K_{1,2}$-isolating set of $M_n$.
Put $U=V(M_n)\setminus N[S]$, so that
$\Delta(M_n[U])\leq 1$.
For $j\in\mathbb Z_m$, define the $j$th
\emph{rung} by $R_j=\{v_j,v_{j+m}\}$.
Rung indices are taken modulo $m$.
A rung is called \emph{empty} if it contains
no vertex of $S$.

\begin{claim}\label{emptyrungs}
   There are at most three consecutive empty rungs.
\end{claim}
\textit{Proof of Claim~\ref{emptyrungs}.}
Suppose, after a cyclic relabelling, that four consecutive rungs
$R_0,R_1,R_2,R_3$ are empty.
Every neighbor of a vertex in $R_1\cup R_2$
belongs to one of these four rungs.
Hence all four vertices of $R_1\cup R_2$
are uncovered.
They induce a $4$-cycle, contradicting
$\Delta(M_n[U])\leq 1$.
The configuration is shown in
Figure~\ref{fig:four-empty-rungs}.
\hfill$\scriptstyle\blacksquare$

\begin{figure}[htbp]
\centering
\begin{tikzpicture}[graph,scale=0.6]
\foreach \i in {0,...,3}{
  \coordinate (t\i) at (1.9*\i,.75);
  \coordinate (b\i) at (1.9*\i,-.75);
  \draw (t\i)--(b\i);
  \node at (1.9*\i,1.35) {$R_{\i}$};
}
\foreach \i in {0,1,2}{
  \pgfmathtruncatemacro{\j}{\i+1}
  \draw (t\i)--(t\j) (b\i)--(b\j);
}
\draw (t0)--++(-.55,0) (b0)--++(-.55,0)
      (t3)--++(.55,0) (b3)--++(.55,0);
\node at (-.9,.75) {$\cdots$};
\node at (-.9,-.75) {$\cdots$};
\node at (6.6,.75) {$\cdots$};
\node at (6.6,-.75) {$\cdots$};

\draw[residualedge] (t1)--(t2)--(b2)--(b1)--cycle;

\foreach \i in {0,3}{
  \node[graphunknown] at (t\i) {};
  \node[graphunknown] at (b\i) {};
}
\foreach \i in {1,2}{
  \node[graphuncovered] at (t\i) {};
  \node[graphuncovered] at (b\i) {};
}

\node[graphuncovered] at (9,0.8) {};
\node[anchor=west] at (9.35,0.8) {uncovered};

\node[graphunknown] at (9,0) {};
\node[anchor=west] at (9.35,0) {coverage unspecified};

\draw[residualedge] (8.75,-0.8)--(9.25,-0.8);
\node[anchor=west] at (9.35,-0.8) {uncovered $C_4$};

\node at (0,-2) {};
\end{tikzpicture}
\caption{Four consecutive empty rungs in $M_n$.}
\label{fig:four-empty-rungs}
\end{figure}
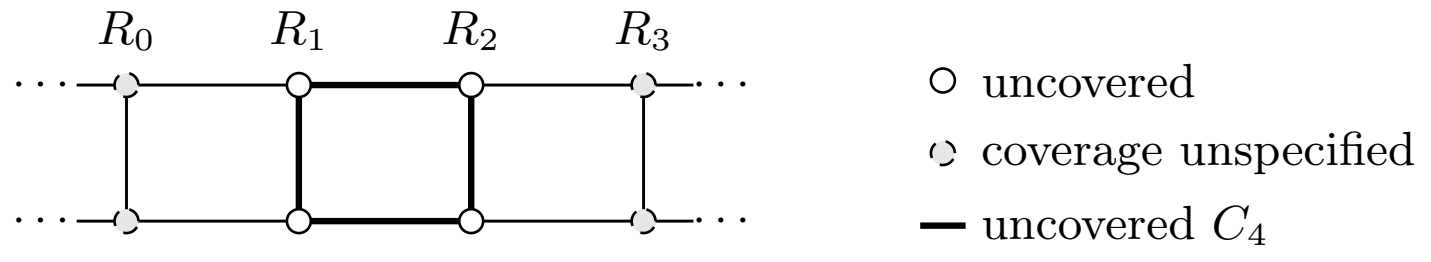

\begin{claim}\label{RS}
If $R_{i-1},R_i,R_{i+1}$ are empty, then
$R_{i-2}\subseteq S$ and $R_{i+2}\subseteq S$.
\end{claim}
\textit{Proof of Claim~\ref{RS}.}
Both vertices of $R_i$ are uncovered and adjacent.
Since $\Delta(M_n[U])\leq 1$, their neighbors
in $R_{i-1}\cup R_{i+1}$ must be covered.
As these rungs are empty, each vertex of
$R_{i-1}$ can only be covered by its neighbor
in $R_{i-2}$, and similarly for $R_{i+1}$.
Thus both bordering rungs belong entirely to $S$.
The conclusion also holds when the bordering
rungs coincide; for $m=3$, three consecutive
empty rungs are impossible.
See Figure~\ref{fig:three-empty-rungs}.
\hfill$\scriptstyle\blacksquare$

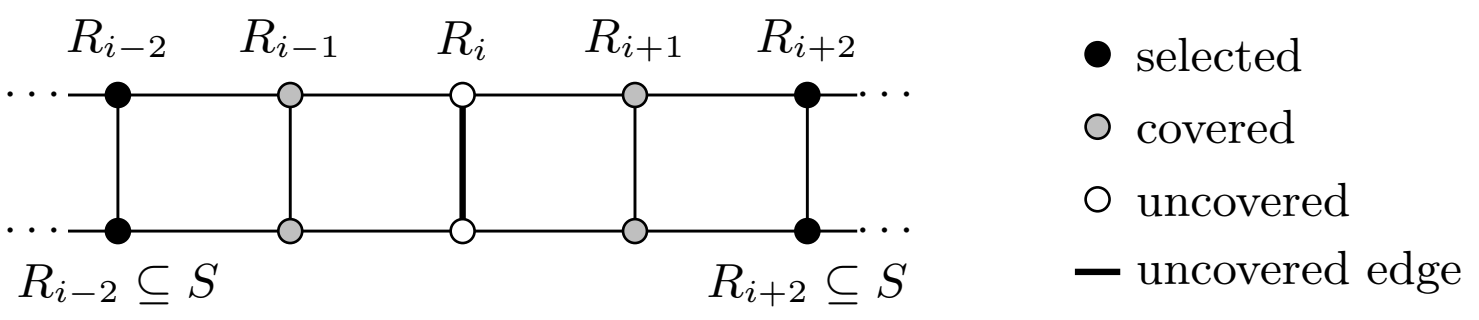
\begin{figure}[htbp]
\centering
\begin{tikzpicture}[graph,scale=0.6]
\foreach \i/\r in {0/{i-2},1/{i-1},2/{i},3/{i+1},4/{i+2}}{
  \coordinate (t\i) at (1.9*\i,.75);
  \coordinate (b\i) at (1.9*\i,-.75);
  \draw (t\i)--(b\i);
  \node at (1.9*\i,1.35) {$R_{\r}$};
}
\foreach \i in {0,1,2,3}{
  \pgfmathtruncatemacro{\j}{\i+1}
  \draw (t\i)--(t\j) (b\i)--(b\j);
}
\draw (t0)--++(-.55,0) (b0)--++(-.55,0)
      (t4)--++(.55,0) (b4)--++(.55,0);

\node at (-.9,.75) {$\cdots$};
\node at (-.9,-.75) {$\cdots$};
\node at (8.5,.75) {$\cdots$};
\node at (8.5,-.75) {$\cdots$};

\draw[residualedge] (t2)--(b2);

\foreach \i in {0,4}{
  \node[graphselected] at (t\i) {};
  \node[graphselected] at (b\i) {};
}
\foreach \i in {1,3}{
  \node[graphcovered] at (t\i) {};
  \node[graphcovered] at (b\i) {};
}
\node[graphuncovered] at (t2) {};
\node[graphuncovered] at (b2) {};

\node at (0,-1.35) {$R_{i-2}\subseteq S$};
\node at (7.6,-1.35) {$R_{i+2}\subseteq S$};

\node[graphselected] at (10.8,1.2) {};
\node[anchor=west] at (11.15,1.2) {selected};

\node[graphcovered] at (10.8,0.4) {};
\node[anchor=west] at (11.15,0.4) {covered};

\node[graphuncovered] at (10.8,-0.4) {};
\node[anchor=west] at (11.15,-0.4) {uncovered};

\draw[residualedge] (10.55,-1.2)--(11.05,-1.2);
\node[anchor=west] at (11.15,-1.2) {uncovered edge};

\node at (0,-2.5) {};

\end{tikzpicture}
\caption{Three empty rungs and their forced bordering selections.}
\label{fig:three-empty-rungs}
\end{figure}

Since $S$ is isolating, it is nonempty.
Partition the cyclic sequence of rungs into blocks,
each consisting of a nonempty rung followed by all
empty rungs up to the next nonempty rung.
Consider a block containing $g$ empty rungs.
By Claim~\ref{emptyrungs}, $g\leq 3$.
If $g\leq 2$, the block contains at most three
rungs and at least one vertex of $S$.
If $g=3$, Claim~\ref{RS} implies that both vertices
of its first rung belong to $S$.
Thus this block contains four rungs and at least
two vertices of $S$.

In either case, the number of rungs in a block
is at most three times the number of its
vertices belonging to $S$.
Summing over all blocks gives $m\leq 3|S|$.
Therefore,
\(
|S|\geq\left\lceil\frac{m}{3}\right\rceil
=\left\lceil\frac{n}{6}\right\rceil.
\)
Since $S$ was arbitrary, this proves the reverse
inequality and completes the proof.
\end{proof}

\section{Equality in the Linear Bounds}
\label{sec:equality}

The bounds in Theorem~\ref{thm:main-bounds}
are both sharp, even within the class of connected
triangle-free cubic graphs.
In this section, equality refers to the linear bounds
$\frac{|V(G)|}{7}$ and $\frac{|V(G)|}{4}$.
We characterize connected cubic graphs satisfying
$\iota_1(G)=\frac{|V(G)|}{7}$ and identify those that are triangle-free.
We then characterize connected bipartite cubic graphs satisfying
$\iota_1(G)=\frac{|V(G)|}{4}$ and establish their vertex-criticality.

\subsection{Equality in the Lower Bound}
\label{sec:lower-bound-equality}

We first examine equality in the lower bound
$\iota_1(G)\geq \frac{|V(G)|}{7}$.
The equality conditions in the counting argument
used to prove Theorem~\ref{thm:main-bounds}
lead to a construction that characterizes all
connected cubic graphs satisfying equality.

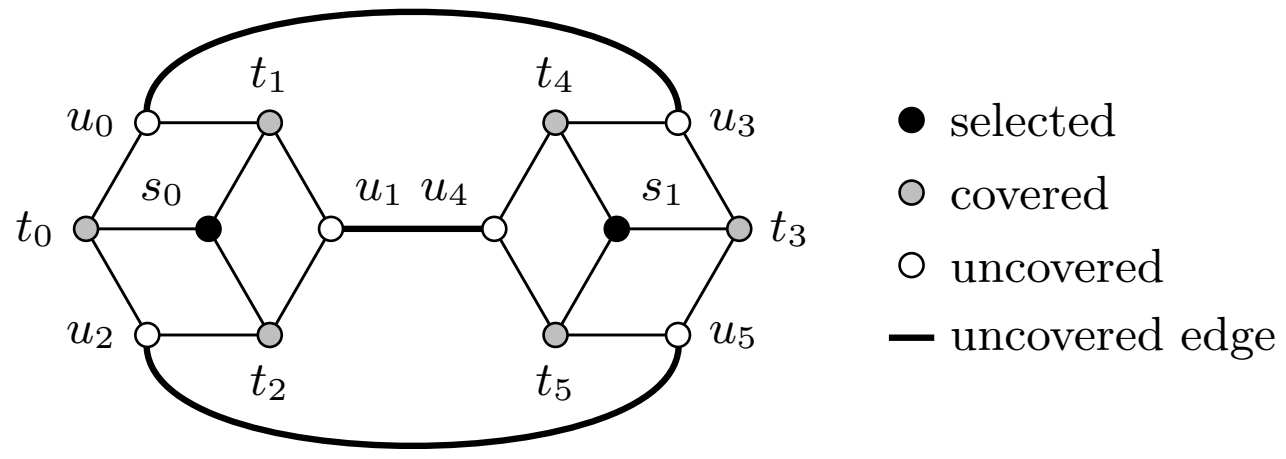
\begin{figure}[htbp]
\centering
\begin{tikzpicture}[graph,scale=0.6]
\coordinate (s0) at (-2.25,0);
\coordinate (s1) at (2.25,0);
\coordinate (t0) at (-3.6,0);
\coordinate (t1) at (-1.575,1.17);
\coordinate (t2) at (-1.575,-1.17);
\coordinate (t3) at (3.6,0);
\coordinate (t4) at (1.575,1.17);
\coordinate (t5) at (1.575,-1.17);
\coordinate (u0) at (-2.925,1.17);
\coordinate (u1) at (-.9,0);
\coordinate (u2) at (-2.925,-1.17);
\coordinate (u3) at (2.925,1.17);
\coordinate (u4) at (.9,0);
\coordinate (u5) at (2.925,-1.17);

\draw (t0)--(u0)--(t1)--(u1)--(t2)--(u2)--cycle;
\draw (t3)--(u3)--(t4)--(u4)--(t5)--(u5)--cycle;

\draw (s0)--(t0) (s0)--(t1) (s0)--(t2)
      (s1)--(t3) (s1)--(t4) (s1)--(t5);

\draw[residualedge]
    (u0) .. controls (-3.2,2.8) and (3.2,2.8) .. (u3);
\draw[residualedge] (u1)--(u4);
\draw[residualedge]
    (u2) .. controls (-3.2,-2.8) and (3.2,-2.8) .. (u5);

\node[graphselected,label=above left:$s_0$] at (s0) {};
\node[graphselected,label=above right:$s_1$] at (s1) {};

\node[graphcovered,label=left:$t_0$] at (t0) {};
\node[graphcovered,label=above:$t_1$] at (t1) {};
\node[graphcovered,label=below:$t_2$] at (t2) {};
\node[graphcovered,label=right:$t_3$] at (t3) {};
\node[graphcovered,label=above:$t_4$] at (t4) {};
\node[graphcovered,label=below:$t_5$] at (t5) {};

\node[graphuncovered,label=left:$u_0$] at (u0) {};
\node[graphuncovered,label=above right:$u_1$] at (u1) {};
\node[graphuncovered,label=left:$u_2$] at (u2) {};
\node[graphuncovered,label=right:$u_3$] at (u3) {};
\node[graphuncovered,label=above left:$u_4$] at (u4) {};
\node[graphuncovered,label=right:$u_5$] at (u5) {};

\node[graphselected] at (5.5,1.2) {};
\node[anchor=west] at (5.85,1.2) {selected};

\node[graphcovered] at (5.5,0.4) {};
\node[anchor=west] at (5.85,0.4) {covered};

\node[graphuncovered] at (5.5,-0.4) {};
\node[anchor=west] at (5.85,-0.4) {uncovered};

\draw[residualedge] (5.25,-1.2)--(5.75,-1.2);
\node[anchor=west] at (5.85,-1.2) {uncovered edge};

\end{tikzpicture}
\caption{A lower-bound equality graph of order fourteen.}
\label{fig:lower-bound-construction}
\end{figure}

For an integer $r\geq 1$, take disjoint sets
$T$ and $U$, each of size $6r$, and let $F$
be a simple $2$-regular bipartite graph with
bipartition $(T,U)$.
Thus $F$ is a disjoint union of alternating
even cycles.
Choose a perfect matching $M$ on $U$.
Partition $T$ into triples
$T_0,\ldots,T_{2r-1}$, introduce new vertices
$S=\{s_0,\ldots,s_{2r-1}\}$, and join each
$s_j$ to all three vertices of $T_j$.
Together with the edges of $F$ and $M$,
these edges form a cubic graph of order $14r$.
We retain only choices for which the resulting graph is connected.
Figure~\ref{fig:lower-bound-construction}
illustrates the construction for $r=1$.

We now show that this construction describes
exactly the connected cubic graphs satisfying $\iota_1(G)=\frac{|V(G)|}{7}$.
The triangle-free members can also be identified
directly from the matching $M$.

\begin{theorem}\label{thm:lower-bound-equality}
Let $G$ be a connected cubic graph of order $n$.
Then $\iota_1(G)=\frac{n}{7}$ if and only if $G$
is obtained by the preceding construction.
In particular, equality requires $14\mid n$.
Moreover, a graph obtained by this construction
is triangle-free if and only if
$N_F(u)\cap N_F(v)=\varnothing$
for every $uv\in M$.
\end{theorem}

\begin{proof}
Suppose that $\iota_1(G)=\frac{n}{7}$, and let $S$
be a minimum $K_{1,2}$-isolating set.
Put $T=N_G[S]\setminus S$ and
$U=V(G)\setminus N_G[S]$, and write
$s=|S|$, $t=|T|$, and $u=|U|$.
As in the proof of Theorem~\ref{thm:main-bounds},
we have $t\leq 3s$ and
$2u\leq e_G(T,U)\leq 2t$.
Consequently,
\(
n=s+t+u\leq s+2t\leq 7s=n,
\)
so $t=u=3s$.
Equality throughout the counting argument
forces every vertex of $S$ to have three
neighbors in $T$, every vertex of $T$
to have exactly one neighbor in $S$
and two in $U$, and every vertex of $U$
to have exactly two neighbors in $T$
and one in $U$.
In particular, both $S$ and $T$ are
independent.

It follows that the edges between $T$
and $U$ form a $2$-regular bipartite graph
$F$, while the edges within $U$ form
a perfect matching $M$.
Moreover, the neighborhoods in $T$
of the vertices of $S$ partition $T$
into triples.
Since $|U|=3s$ is even, $s=2r$
for some integer $r\geq 1$.
Thus $G$ has the stated construction.

Conversely, suppose that $G$ is obtained
by this construction.
The set $S$ is $K_{1,2}$-isolating,
since $G-N_G[S]$ consists precisely
of the matching $M$.
Hence $\iota_1(G)\leq |S|=2r=\frac{n}{7}$.
The counting argument in the proof of Theorem~\ref{thm:main-bounds},
which applies to every connected cubic graph,
gives the reverse inequality, proving equality.

Finally, no triangle contains a vertex of $S$,
since every neighbor of a vertex of $S$ lies in the independent set $T$.
As $T$ is also independent, every triangle
must consist of an edge $uv\in M$
and a common neighbor of $u$ and $v$ in $T$.
Thus $G$ is triangle-free precisely when
$N_F(u)\cap N_F(v)=\varnothing$
for every $uv\in M$.
\end{proof}
\noindent
The construction includes connected triangle-free
cubic graphs, as shown in
Figure~\ref{fig:lower-bound-construction}.
Hence the lower bound is sharp even within
this class.

\medskip

Furthermore,
the Coxeter graph provides another example.
This well-known cubic graph, described by
Tutte~\cite{tutte1960non}, has order $28$
and girth seven.

\begin{corollary}\label{cor:coxeter-lower-equality}
The Coxeter graph $CG$ belongs to the construction
of Theorem~\ref{thm:lower-bound-equality}
with $r=2$.
In particular, $\iota_1(CG)=4$.
\end{corollary}
\noindent
The partition $(S,T,U)$ is an equitable partition
(a perfect coloring), and $S$ is a completely regular code
of covering radius two~\cite{DamascoFrettloh2021,KrotovPotapov2025}.
These standard notions give an alternative interpretation
of the equality case $\iota_1(G)=\frac{n}{7}$.

\subsection{Bipartite Graphs Attaining the Upper Bound}
\label{sec:bipartite-equality}

We now investigate equality in the upper bound
$\iota_1(G)\leq \frac{|V(G)|}{4}$ for connected bipartite cubic graphs.
The bound follows from the independent-isolation result of
Boyer and Goddard~\cite{BoyerGoddard2025};
our focus here is its equality cases.
Our approach builds on known results concerning
total domination, particularly an infinite family
of cubic graphs studied by Favaron et
al.~\cite{FavaronEtAl2000}, Henning and
Yeo~\cite{HenningYeo2008}, and
Cyman et al.~\cite{CymanEtAl2018}.
We begin by showing that this family attains
equality for $K_{1,2}$-isolation.

\medskip

\begin{figure}[htbp]
\centering
\begin{tikzpicture}[graph,scale=0.6]
\foreach \i in {0,1,2}{
  \coordinate (a\i) at (3.1*\i,.75);
  \coordinate (b\i) at (3.1*\i+1.4,.75);
  \coordinate (c\i) at (3.1*\i,-.75);
  \coordinate (d\i) at (3.1*\i+1.4,-.75);
  \draw (a\i)--(b\i)--(c\i)--(d\i)--cycle;
  \node[graphselected] at (a\i) {};
  \node[graphcovered] at (b\i) {};
  \node[graphuncovered] at (c\i) {};
  \node[graphcovered] at (d\i) {};
}
\draw (b0)--(a1) (b1)--(a2) (d0)--(c1) (d1)--(c2);
\draw (b2) .. controls (8.2,2.8) and (-.6,2.8) .. (a0);
\draw (d2) .. controls (8.2,-2.8) and (-.6,-2.8) .. (c0);

\node[left=4pt] at (a0) {$a_0$};
\node[left=4pt] at (c0) {$c_0$};
\node[right=4pt] at (b2) {$b_2$};
\node[right=4pt] at (d2) {$d_2$};

\foreach \i in {1,2}{
  \node[above=4pt] at (a\i) {$a_{\i}$};
  \node[below=4pt] at (c\i) {$c_{\i}$};
}
\foreach \i in {0,1}{
  \node[above=4pt] at (b\i) {$b_{\i}$};
  \node[below=4pt] at (d\i) {$d_{\i}$};
}

\node[graphselected] at (10,1.2) {};
\node[anchor=west] at (10.35,1.2) {selected};

\node[graphcovered] at (10,0.4) {};
\node[anchor=west] at (10.35,0.4) {covered};

\node[graphuncovered] at (10,-0.4) {};
\node[anchor=west] at (10.35,-0.4) {uncovered};

\draw[residualedge] (9.75,-1.2)--(10.25,-1.2);
\node[anchor=west] at (10.35,-1.2) {uncovered edge};

\end{tikzpicture}
\caption{$H_3$ with $S=\{a_0,a_1,a_2\}$.}
\label{fig:infinite-equality}
\end{figure}
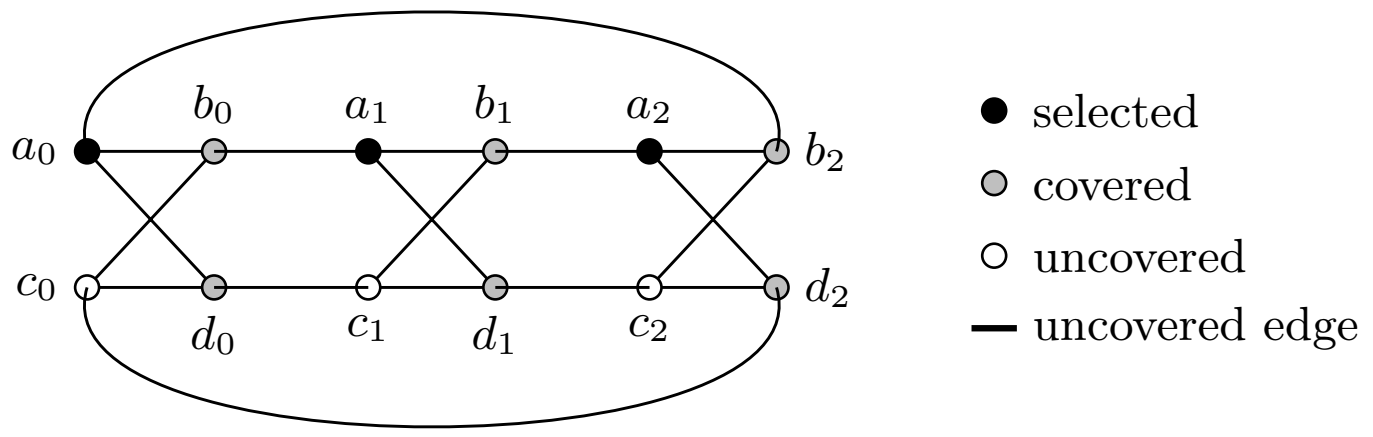

For $k\geq2$, let $H_k$ be the graph obtained
from $k$ vertex-disjoint four-cycles
\(
Q_i=a_i b_i c_i d_i a_i,
\quad i=0,\ldots,k-1,
\)
by adding the edges $b_i a_{i+1}$ and
$d_i c_{i+1}$, with indices taken modulo $k$.
Let $\mathcal H=\{H_k:k\geq2\}$.
Put $U_i=\{a_i,c_i\}$ and $V_i=\{b_i,d_i\}$.
Then $H_k$ is a connected bipartite cubic graph
with bipartition
\(
\left(\bigcup_{i=0}^{k-1}U_i,\,
      \bigcup_{i=0}^{k-1}V_i\right).
\)
Moreover, $|V(H_k)|=4k$ and $H_2\cong Q_3$, the cube graph.
Our notation agrees with the family
$H_\ell$ in~\cite{CymanEtAl2018},
with $\ell=k$.
Figure~\ref{fig:infinite-equality}
illustrates $H_3$.

\medskip

We first establish that every graph in
$\mathcal H$ attains the upper bound.

\begin{theorem}\label{thm:infinite-equality}
For every $k\geq2$, $\iota_1(H_k)=k=\frac{|V(H_k)|}{4}$.
\end{theorem}

\begin{proof}
Since $H_k$ is a triangle-free cubic graph of order $4k$,
Theorem~\ref{thm:main-bounds} gives
$\iota_1(H_k)\leq k$.

For the reverse inequality, let $S$ be any $K_{1,2}$-isolating set.
Write $x_i=|S\cap U_i|$, $y_i=|S\cap V_i|$, and $s_i=x_i+y_i$.
Call $Q_i$ occupied if $s_i\geq1$, and empty otherwise.
In an empty $Q_i$, exactly $y_{i-1}$ vertices of $U_i$ and
$x_{i+1}$ vertices of $V_i$ are covered: the external attachment
edges are matchings. At least two vertices of this four-cycle
must be covered, since deleting at most one vertex leaves a $P_3$.
Consequently,
\(
y_{i-1}+x_{i+1}\geq2
\text{ whenever }s_i=0.
\)

Give each empty cycle one unit of charge, and distribute it to
its neighboring occupied cycles as follows. If $y_{i-1}=2$, give
all the charge to $Q_{i-1}$. Otherwise, if $x_{i+1}=2$, give it
all to $Q_{i+1}$. In the remaining case,
$y_{i-1}=x_{i+1}=1$; give half to each neighbor.

For $k\geq3$, consider the split case
$y_{i-1}=x_{i+1}=1$.
The empty block $Q_i$ has exactly two uncovered vertices,
say $u\in U_i$ and $w\in V_i$, which are adjacent.
Let $z\in V_{i-1}$ be the external neighbor of $u$.
Since $u$ is uncovered, $z$ is unselected.
If $x_{i-1}=0$, neither of the two neighbors of $z$ in
$U_{i-1}$ is selected, and its remaining neighbor $u$
lies in the empty block.
Thus $z$ is also uncovered, leaving an uncovered $P_3$ on $z,u,w$,
a contradiction.
Hence $x_{i-1}\geq1$; the symmetric argument gives $y_{i+1}\geq1$.
Therefore, each block receiving half a unit contains selections
in both parts.
For $k=2$, the two neighboring blocks coincide, and the split
condition itself gives one selection in each part of that block.
Figure~\ref{fig:half-charge-obstruction} illustrates the obstruction.

\begin{figure}[htbp]
\centering
\begin{tikzpicture}[graph,scale=0.8]
\foreach \j/\r in {0/{i-1},1/{i},2/{i+1}}{
  \coordinate (a\j) at (3.5*\j,.8);
  \coordinate (b\j) at (3.5*\j+1.4,.8);
  \coordinate (c\j) at (3.5*\j,-.8);
  \coordinate (d\j) at (3.5*\j+1.4,-.8);
  \draw[graphguide] (3.5*\j-.55,-1.35) rectangle (3.5*\j+1.95,1.35);
  \draw (a\j)--(b\j)--(c\j)--(d\j)--cycle;
  \node at (3.5*\j+.7,1.85) {$Q_{\r}$};
}

\draw (b0)--(a1) (b1)--(a2) (d0)--(c1) (d1)--(c2);
\draw[residualedge] (d0)--(c1)--(d1);

\node[graphcovered,label=above:$a_{i-1}$] at (a0) {};
\node[graphselected,label=above:$b_{i-1}$] at (b0) {};
\node[graphcovered,label=below:$c_{i-1}$] at (c0) {};
\node[graphuncovered,label=below:$d_{i-1}$] at (d0) {};

\node[graphcovered,label=above:$a_i$] at (a1) {};
\node[graphcovered,label=above:$b_i$] at (b1) {};
\node[graphuncovered,label=below:$c_i$] at (c1) {};
\node[graphuncovered,label=below:$d_i$] at (d1) {};

\node[graphselected,label=above:$a_{i+1}$] at (a2) {};
\node[graphcovered,label=above:$b_{i+1}$] at (b2) {};
\node[graphuncovered,label=below:$c_{i+1}$] at (c2) {};
\node[graphcovered,label=below:$d_{i+1}$] at (d2) {};

\node at (.7,-1.85) {$x_{i-1}=0,\ y_{i-1}=1$};
\node at (4.2,-1.85) {$s_i=0$};
\node at (7.7,-1.85) {$x_{i+1}=1$};

\node[graphselected] at (10.2,1.2) {};
\node[anchor=west] at (10.55,1.2) {selected};

\node[graphcovered] at (10.2,0.4) {};
\node[anchor=west] at (10.55,0.4) {covered};

\node[graphuncovered] at (10.2,-0.4) {};
\node[anchor=west] at (10.55,-0.4) {uncovered};

\draw[residualedge] (9.95,-1.2)--(10.45,-1.2);
\node[anchor=west] at (10.55,-1.2) {uncovered $P_3$};

\end{tikzpicture}
\caption{An uncovered $P_3$ when $x_{i-1}=0$ in the split case.}
\label{fig:half-charge-obstruction}
\end{figure}
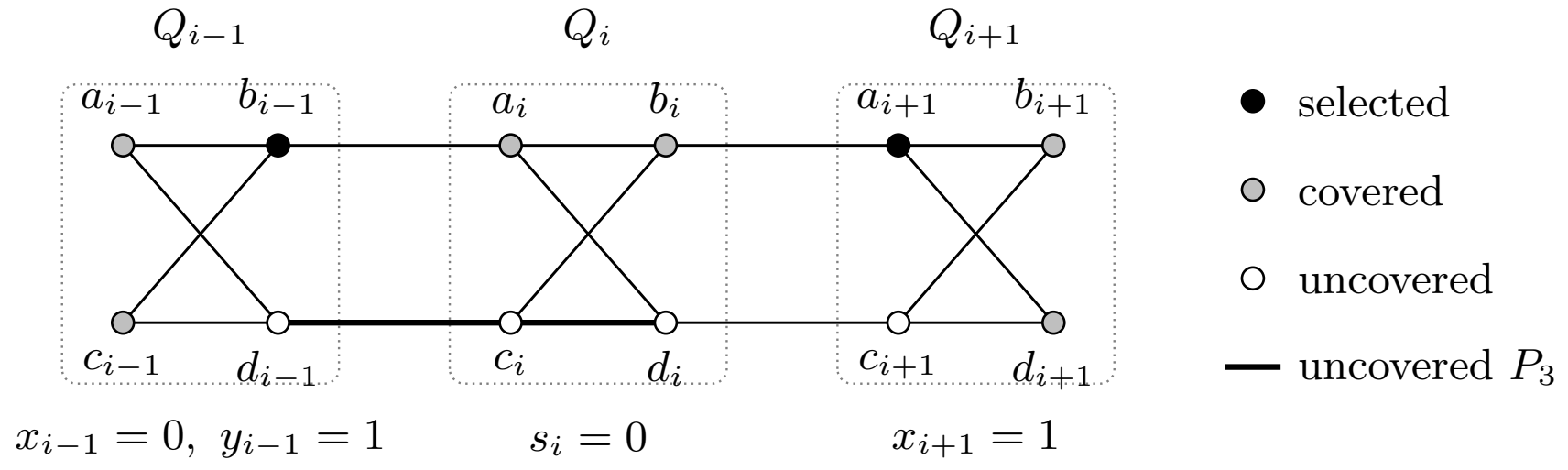

We claim that each occupied $Q_j$ receives at most $s_j-1$ units.
If $s_j=1$, it receives neither a full charge, which requires two
selections in the same part, nor a half-charge, which requires
selections in both parts. If $s_j=2$ and both selections lie in
the same part, it can receive a full charge only from the one
corresponding direction, and cannot receive a half-charge.
If $s_j=2$ with one selection in each part, it can receive at
most one half-charge from each direction, and no full charge.
Finally, if $s_j\geq3$, it receives at most one unit from each
direction, hence at most $2\leq s_j-1$ units.
For $k=2$ the two directions lead to the same neighboring cycle;
the same bounds apply, counting the two possible half-charges
separately.

Let $e$ and $o$ be the numbers of empty and occupied cycles.
Summing the charges gives
\[
e\leq\sum_{j:s_j>0}(s_j-1)=|S|-o.
\]
Since $k=e+o$, we obtain $|S|\geq k$, as required.
\end{proof}

Having established equality for the infinite family
$\mathcal H$, we next determine whether other connected
bipartite cubic graphs attain the upper bound.
For this purpose, we use a characterization of cubic
graphs with total domination number equal to half
their order, due to Henning and Yeo~\cite{HenningYeo2008}.
In the bipartite case, their characterization identifies
one additional graph outside $\mathcal H$:
the M\"obius-Kantor graph.

\medskip

\begin{figure}[htbp]
\centering
\begin{tikzpicture}[graph,scale=0.7]
\foreach \i in {0,...,7}{
  \pgfmathsetmacro{\ang}{90-45*\i}
  \pgfmathsetmacro{\labelang}{\ang+14}
  \coordinate (u\i) at (\ang:2.25);
  \coordinate (w\i) at (\ang:1.15);
  \coordinate (ul\i) at (\ang:2.75);
  \coordinate (wl\i) at (\labelang:1.6);
}
\foreach \i in {0,...,7}{
  \pgfmathtruncatemacro{\j}{mod(\i+1,8)}
  \pgfmathtruncatemacro{\k}{mod(\i+3,8)}
  \draw (u\i)--(u\j) (u\i)--(w\i) (w\i)--(w\k);
}
\foreach \i in {0,2,4,6}{
  \node[graphselected] at (u\i) {};
  \node[graphcovered] at (w\i) {};
}
\foreach \i in {1,3,5,7}{
  \node[graphcovered] at (u\i) {};
  \node[graphuncovered] at (w\i) {};
}
\foreach \i in {0,...,7}{
  \node at (ul\i) {$u_{\i}$};
  \node at (wl\i) {$w_{\i}$};
}

\node[graphselected] at (4.5,0.8) {};
\node[anchor=west] at (4.85,0.8) {selected};

\node[graphcovered] at (4.5,0) {};
\node[anchor=west] at (4.85,0) {covered};

\node[graphuncovered] at (4.5,-0.8) {};
\node[anchor=west] at (4.85,-0.8) {uncovered};

\end{tikzpicture}
\caption{A minimum $K_{1,2}$-isolating set of the M\"obius--Kantor graph.}
\label{fig:mobius-kantor-witness}
\end{figure}
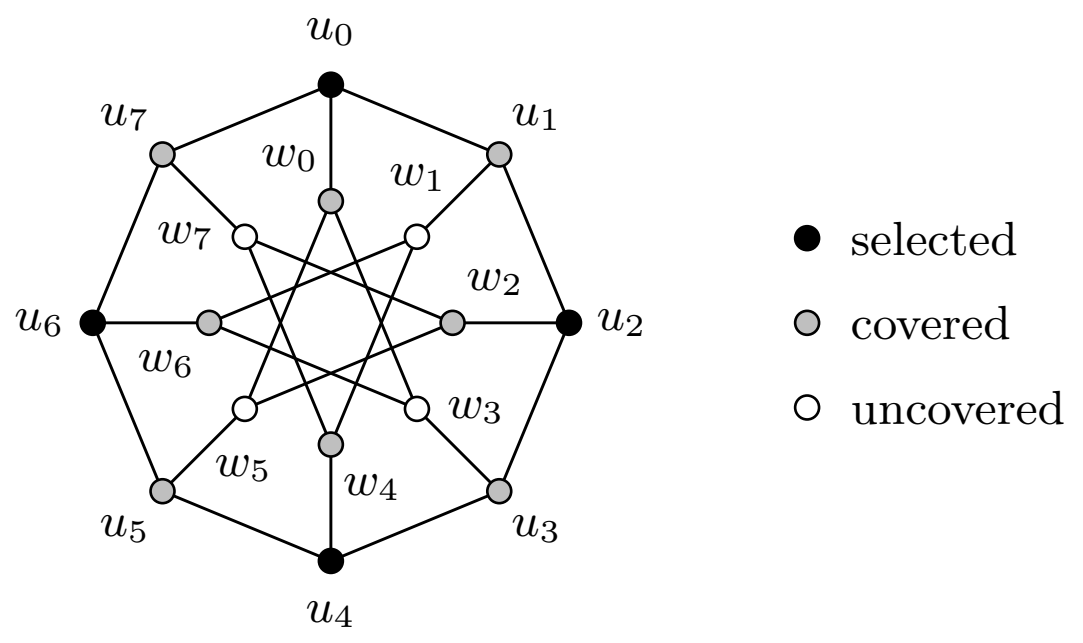

The M\"obius-Kantor graph
$M=\operatorname{GP}(8,3)$ is a connected
bipartite cubic graph of order sixteen,
shown in Figure~\ref{fig:mobius-kantor-witness}.
Its vertex set is
$\{u_i,w_i:i\in\mathbb Z_8\}$, with edges
$u_iu_{i+1}$, $u_iw_i$, and $w_iw_{i+3}$,
where indices are taken modulo eight.
A bipartition of $M$ is given by
\(
A=\{u_i:i\text{ even}\}
   \cup\{w_i:i\text{ odd}\}
\)
and
\(
B=\{u_i:i\text{ odd}\}
   \cup\{w_i:i\text{ even}\}.
\)

The relevant bipartite case of their characterization
can be stated as follows.

\begin{proposition}[\cite{HenningYeo2008}]
\label{prop:bipartite-td-equality}
Let $G$ be a connected bipartite cubic graph
of order $n$.
Then $\gamma_{\mathrm t}(G)=\frac{n}{2}$ if and only if
$G\cong H_k$ for some $k\geq2$ or
$G\cong M$.
\end{proposition}

We next verify that the M\"obius-Kantor graph
also attains equality for $K_{1,2}$-isolation.

\begin{proposition}
\label{prop:mobius-kantor-isolation}
The M\"obius-Kantor graph satisfies
$\iota_1(M)=4$.
\end{proposition}

\begin{proof}
The edges induced by the $u$-vertices form
an eight-cycle, as do those induced by the
$w$-vertices, since $\gcd(8,3)=1$.
The spokes $u_iw_i$ form a matching.
Any four-cycle containing spokes would
therefore have exactly two of them,
together with one edge from each
eight-cycle.
This would require two indices to differ
by both $\pm1$ and $\pm3$ modulo eight,
which is impossible.
Hence $M$ contains no four-cycle.

Let $A,B$ be the bipartition of $M$,
with $|A|=|B|=8$.
Since $M$ is cubic and contains no
four-cycle, any two vertices in the same
part have at most one common neighbor.
Moreover, for every $p\in A$, its three
neighbors have six distinct neighbors
in $A\setminus\{p\}$.
Thus exactly one vertex of
$A\setminus\{p\}$ has no common neighbor
with $p$.
Consequently, the pairs of vertices
with disjoint neighborhoods form a
perfect matching on $A$, and similarly
on $B$.

Suppose that $S$ is a $K_{1,2}$-isolating
set of cardinality three.
Interchanging $A$ and $B$ if necessary,
we may assume that either
$S\subseteq A$ or
$|S\cap A|=2$ and $|S\cap B|=1$.
First suppose $S=\{p,q,r\}\subseteq A$.
Since the pairs with disjoint
neighborhoods form a matching, at least
two of the three pairs among $p,q,r$
have a common neighbor.
Writing
\(
h=
|N(p)\cap N(q)|
+|N(p)\cap N(r)|
+|N(q)\cap N(r)|
\)
and
\(
t=|N(p)\cap N(q)\cap N(r)|,
\)
we have $h\geq2$.
If $t=0$, then $h-t\geq2$; if $t=1$,
all three pairwise intersections have
size one, so $h=3$ and again $h-t\geq2$.
Inclusion--exclusion therefore gives
$|N(S)|=9-h+t\leq7$.
Hence some vertex $v\in B$ is uncovered.
Since no vertex of $A\setminus S$
is covered, all three neighbors of $v$
are also uncovered, contradicting
the assumption that $S$ is
$K_{1,2}$-isolating.

Now suppose $S=\{p,q,z\}$,
where $p,q\in A$ and $z\in B$.
At most seven vertices of $B$
belong to $N[S]$, since
$N[S]\cap B\subseteq
N(p)\cup N(q)\cup\{z\}$.
Choose an uncovered vertex $v\in B$.
Since $v$ and $z$ have at most one
common neighbor, at least two neighbors
of $v$ lie outside $N(z)$.
Neither can belong to $S$, because
$v$ has no selected neighbor.
Thus these vertices and $v$ form
an uncovered $P_3$, again a contradiction.

Therefore, no set of cardinality three
is $K_{1,2}$-isolating, and
$\iota_1(M)\geq4$.
On the other hand, since $M$ is a
triangle-free cubic graph of order sixteen,
Theorem~\ref{thm:main-bounds}
gives $\iota_1(M)\leq4$.
Consequently, $\iota_1(M)=4$.
\end{proof}

The total domination characterization in
Proposition~\ref{prop:bipartite-td-equality} gives the necessary candidates
for $\iota_1(G)=\frac{|V(G)|}{4}$.
The preceding $K_{1,2}$-isolation lower bounds establish sufficiency,
yielding the following characterization.

\begin{theorem}\label{thm:bipartite-equality}
Let $G$ be a connected bipartite cubic
graph of order $n$.
Then $\iota_1(G)=\frac{n}{4}$ if and only if
$G\cong H_k$ for some $k\geq2$
or $G\cong M$.
\end{theorem}

\begin{proof}
Let $A,B$ be the bipartition, and let
$D$ be a minimum total dominating set
of $G$.
Choose the smaller of $D\cap A$
and $D\cap B$, say $S=D\cap A$.
Every vertex of $B$ has a neighbor
in $S$, since $D$ is total dominating.
Thus every uncovered vertex lies
in the independent part $A$, so
$S$ is even an ordinary isolating set.
By Proposition~
\ref{prop:td-short-degree-two-components},
\[
\iota_1(G)
\leq |S|
\leq\frac{\gamma_{\mathrm t}(G)}2
\leq\frac n4.
\]
If $\iota_1(G)=\frac{n}{4}$, the above
inequalities imply that
$\gamma_{\mathrm t}(G)=\frac{n}{2}$.
Proposition~\ref{prop:bipartite-td-equality}
now gives the required list.
Conversely, every graph on that list
attains equality by
Theorem~\ref{thm:infinite-equality}
and
Proposition~\ref{prop:mobius-kantor-isolation}.
\end{proof}

We conclude by examining the effect of vertex deletion
on the $K_{1,2}$-isolation number.
The following result gives a general bound for
connected bipartite cubic graphs and, in particular, shows
that every such graph satisfying $\iota_1(G)=\frac{|V(G)|}{4}$
is vertex-critical.

\begin{proposition}
\label{prop:bipartite-vertex-deletion}
Let $G$ be a connected bipartite cubic graph
of order $n$.
For every $v\in V(G)$,
\[
\iota_1(G-v)
\leq
\left\lfloor\frac{n-1}{4}\right\rfloor.
\]
If $\iota_1(G)=\frac{n}{4}$, then
$\iota_1(G-v)=\iota_1(G)-1$
for every $v\in V(G)$.
\end{proposition}
\begin{proof}
Put $F=G-v$.
Its vertices have degree two or three,
and its degree-two vertices are precisely
$N_G(v)$, an independent set.
Apply Proposition~\ref{prop:td-short-degree-two-components}
separately to each component $C$ of $F$, obtaining a total
dominating set $D_C$ with $|D_C|\leq \frac{|V(C)|}{2}$.
The union $D$ of these sets totally dominates $F$ and satisfies
$|D|\leq\frac{n-1}{2}$.
For a bipartition $A,B$ of $F$,
choose the smaller of $D\cap A$
and $D\cap B$, say $X=D\cap A$.
Every vertex of $B$ has a neighbor
in $X$, so the uncovered vertices
lie in the independent set $A$.
Therefore,
\[
\iota_1(G-v)
\leq |X|
\leq\frac{|D|}{2}
\leq\frac{n-1}{4},
\]
and integrality gives the first assertion.
Conversely, adjoining $v$ to a
$K_{1,2}$-isolating set of $G-v$ gives one
for $G$, so
$\iota_1(G-v)\geq\iota_1(G)-1$.
If $\iota_1(G)=\frac{n}{4}$,
the two bounds coincide.
\end{proof}
\section{Concluding Remarks}
\label{sec:conclusion}

We have established sharp lower and upper bounds
for the $K_{1,2}$-isolation number of connected triangle-free
cubic graphs. The upper bound confirms the triangle-free
case of Conjecture~\ref{conj:cubic-path-isolation},
while the lower bound holds more generally for all connected
cubic graphs. We have characterized the connected cubic graphs
satisfying $\iota_1(G)=\frac{|V(G)|}{7}$ and the connected bipartite
cubic graphs satisfying $\iota_1(G)=\frac{|V(G)|}{4}$, showing
that the latter are vertex-critical.

\medskip

Our upper bound, together with the results of
Yin, An and Wu~\cite{YinAnWu2024} and
Bartolo, Borg and Scicluna~\cite{BartoloBorgScicluna2026},
further narrows the possible counterexamples
to Conjecture~\ref{conj:cubic-path-isolation}.
Any such counterexample must contain a triangle,
an induced claw, and an induced six-cycle.
Although the connected bipartite cubic graphs satisfying
$\iota_1(G)=\frac{|V(G)|}{4}$ have been completely characterized,
this equality is not confined to bipartite graphs.
We demonstrate this by constructing a connected
nonbipartite triangle-free cubic graph of order
sixteen and girth five that also attains equality.

\medskip

Let $J$ be obtained from two disjoint eight-cycles
$x_0x_1\cdots x_7x_0$ and
$y_0y_1\cdots y_7y_0$
by adding the chords
$x_1x_5$, $x_3x_7$, $y_1y_5$, and $y_3y_7$,
together with the edges $x_iy_i$
for $i\in\{0,2,4,6\}$.
Equivalently, $J$ can be viewed as joining
two copies of the underlying graph of the
Petersen $4$-pole considered by
H\"agglund~\cite{Hagglund2016}.

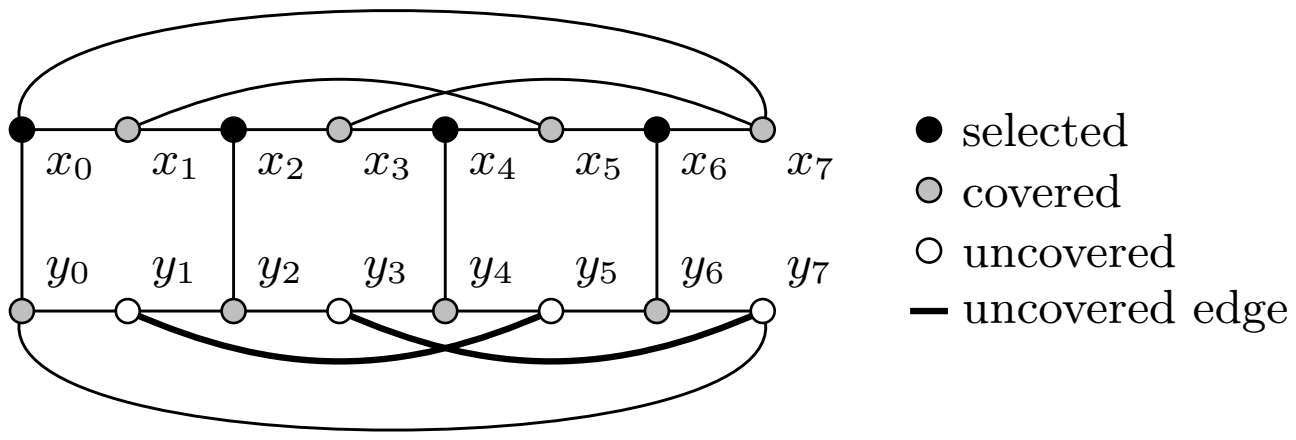
\begin{figure}[htbp]
\centering
\begin{tikzpicture}[graph,scale=0.5]
\foreach \i in {0,...,7}{
  \coordinate (x\i) at (1.4*\i,1.2);
  \coordinate (y\i) at (1.4*\i,-1.2);
}
\foreach \i in {0,...,6}{
  \pgfmathtruncatemacro{\j}{\i+1}
  \draw (x\i)--(x\j) (y\i)--(y\j);
}
\draw (x7) .. controls (10.6,3.3) and (-.8,3.3) .. (x0);
\draw (y7) .. controls (10.6,-3.3) and (-.8,-3.3) .. (y0);
\draw (x1) to[bend left=24] (x5);
\draw (x3) to[bend left=24] (x7);
\draw[residualedge] (y1) to[bend right=24] (y5);
\draw[residualedge] (y3) to[bend right=24] (y7);

\foreach \i in {0,2,4,6}{
  \draw (x\i)--(y\i);
  \node[graphselected,label=below right:$x_{\i}$] at (x\i) {};
  \node[graphcovered,label=above right:$y_{\i}$] at (y\i) {};
}
\foreach \i in {1,3,5,7}{
  \node[graphcovered,label=below right:$x_{\i}$] at (x\i) {};
  \node[graphuncovered,label=above right:$y_{\i}$] at (y\i) {};
}

\node[graphselected] at (12,1.2) {};
\node[anchor=west] at (12.35,1.2) {selected};

\node[graphcovered] at (12,0.4) {};
\node[anchor=west] at (12.35,0.4) {covered};

\node[graphuncovered] at (12,-0.4) {};
\node[anchor=west] at (12.35,-0.4) {uncovered};

\draw[residualedge] (11.75,-1.2)--(12.25,-1.2);
\node[anchor=west] at (12.35,-1.2) {uncovered edge};

\end{tikzpicture}
\caption{A minimum $K_{1,2}$-isolating set of $J$.}
\label{fig:nonbipartite-J}
\end{figure}

\begin{proposition}\label{prop:nonbipartite-J}
The graph $J$ is a connected nonbipartite cubic graph
of order sixteen and girth five, and
$
\iota_1(J)=4=\frac{|V(J)|}{4}.
$
\end{proposition}
\begin{proof}
The graph $J$ has order sixteen and is connected and cubic.
Within either eight-vertex copy, a cycle using one chord
has length at least five, while a cycle using both chords
has length at least six. Any cycle meeting both copies
uses at least two joining edges, whose endpoints are nonadjacent
within each copy, and therefore has length at least six.
Since $x_1x_2x_3x_4x_5x_1$ is a five-cycle,
$J$ has girth five and is nonbipartite.
The set $S=\{x_0,x_2,x_4,x_6\}$ is
$K_{1,2}$-isolating, since the only edges
remaining outside $N_J[S]$ are
$y_1y_5$ and $y_3y_7$.
Thus $\iota_1(J)\leq4$.

To prove the reverse inequality, suppose that
$S'$ is a $K_{1,2}$-isolating set of size three.
If all three vertices belong to the same eight-cycle,
say the $x$-cycle, at most three of $y_0,y_2,y_4,y_6$
are covered, while every odd-indexed $y$-vertex is uncovered.
An uncovered even-indexed $y$-vertex and its two cycle neighbors
therefore form an uncovered $P_3$.
Otherwise, using the automorphisms obtained
by interchanging the two cycles and simultaneously
shifting all indices by an even number,
we may assume that
$S'=\{y_j,x_p,x_q\}$ with $j\in\{0,1\}$.
If $j=1$, the vertices of the $y$-cycle
outside $N_J[y_1]$ induce the path
$y_4y_3y_7y_6$.
Its internal vertices $y_3$ and $y_7$ cannot be covered
by either selected $x$-vertex.
Both endpoints must therefore be covered, forcing
$\{x_p,x_q\}=\{x_4,x_6\}$.
But then $x_0x_1x_2$ remains uncovered,
a contradiction.
If $j=0$, the vertices $y_3$ and $y_5$ remain uncovered.
To avoid the uncovered path $y_3y_4y_5$, the vertex $y_4$
must be covered, forcing $x_4\in S'$.
Write $S'=\{y_0,x_4,x_r\}$.
The path $x_2y_2y_3$ remains uncovered unless
$r\in\{1,2,3\}$, while $x_6y_6y_5$ remains uncovered
unless $r\in\{5,6,7\}$.
These requirements are incompatible.
Hence no three-vertex set is
$K_{1,2}$-isolating, and therefore
$\iota_1(J)=4=|V(J)|/4$.
\end{proof}
\noindent
The graph and a minimum $K_{1,2}$-isolating set
are illustrated in Figure~\ref{fig:nonbipartite-J}.

\noindent
This example shows that the upper bound can
be attained by nonbipartite triangle-free
cubic graphs, raising the question of whether
these equality graphs admit a complete
structural characterization.

\begin{problem}\label{prob:nonbipartite-equality}
Characterize all connected nonbipartite
triangle-free cubic graphs $G$ satisfying
$\iota_1(G)=\frac{|V(G)|}{4}$.
\end{problem}

A further question concerns vertex-criticality.
Proposition~\ref{prop:bipartite-vertex-deletion} shows that every connected bipartite
cubic graph satisfying $\iota_1(G)=\frac{|V(G)|}{4}$ loses
exactly one from its $K_{1,2}$-isolation number
upon deletion of any vertex.
It remains to determine whether the same
property holds without the bipartiteness
assumption.

\begin{problem}\label{prob:equality-criticality}
Is every connected triangle-free cubic graph $G$ satisfying
$\iota_1(G)=\frac{|V(G)|}{4}$ vertex-critical; that is, does
$\iota_1(G-v)=\iota_1(G)-1$ hold for every $v\in V(G)$?
\end{problem}

\bibliographystyle{elsarticle-num}
\bibliography{references}
\end{document}